# The Poisson Transmuted-G Family of Distributions: Its Properties and Applications

[1] Laba Handique, [2] Subrata Chakraborty

[1,2] Department of Statistics, Dibrugarh University, Dibrugarh-786004, India

## Abstract

In this paper introduces a new family of continuous distributions namely the Poison transmuted-G family of distribution is proposed by inducing two addition parameter on the base line G distribution. Some of its mathematical properties including explicit expressions for the moments generating function, order statistics, Probability weighted moments, stress-strength reliability, residual life, reversed residual life, Rényi entropy and mean deviation are derived. Some special models of the new family are listed. Estimation of the model parameters by the maximum likelihood method is discussed. The advantage of the proposed family in data fitting is illustrated by means of two applications to failure time data set.

## 1. Introduction

Recently, several generalized families of continuous distributions have been proposed by extending a classical probability distribution and applied to model various phenomena. However, there is a clear need for extended forms of the well-known distributions by adding one or more shape parameter(s) in order to obtain greater flexibility in modelling various data.

Shaw and Buckley (2007) proposed the transmuted-G (T-G) family with cdf and pdf respectively given by

$$F^{TG}(x;\alpha) = G(x)\,[1+\alpha - \alpha\,G(x)] \qquad (1)$$

and
$$f^{TG}(x;\alpha) = g(x)\,[1+\alpha - 2\alpha\,G(x)] \qquad (2)$$

where $G(x)$ and $g(x)$ are the baseline cdf and pdf respectively. For $\alpha = 0$, $(|\alpha| \leq 1)$ then Eq. (1) gives the baseline distribution.

The Poisson-G family of distribution with cdf is given by (see Kumaraswamy Poisson-G family, Chakraborty et al., 2020)

$$F^{PG}(x;\beta) = \frac{1-e^{-\beta F(x)}}{1-e^{-\beta}}, \; \beta \in R-\{0\}\,;\; n=1, 2,... \qquad (3)$$

The corresponding pdf of the Poisson-G family is given by

$$f^{PG}(x;\beta) = \frac{\beta\,f(x)\,e^{-\beta F(x)}}{1-e^{-\beta}}, \; \beta \in R-\{0\}\,;\; -\infty < x < \infty \qquad (4)$$

---

In this article, we introduce a new extension of T-G distribution having two parameters $\alpha$ and $\beta$ by considering the T-G as the baseline distribution in the Poisson-G family of distribution. The pdf, cdf and hrf of this proposed distribution are respectively given by:

$$f^{\text{PTG}}(x;\alpha,\beta) = \frac{\beta g(x)[1+\alpha-2\alpha G(x)]\exp[-\beta G(x)\{1+\alpha-\alpha G(x)\}]}{1-e^{-\beta}} \quad (5)$$

$$F^{\text{PTG}}(x;\alpha,\beta) = \frac{1-\exp[-\beta G(x)\{1+\alpha-\alpha G(x)\}]}{1-e^{-\beta}} \quad (6)$$

and 

$$h^{\text{PTG}}(x;\alpha,\beta) = \frac{\beta g(x)[1+\alpha-2\alpha G(x)]\exp[-\beta G(x)\{1+\alpha-\alpha G(x)\}]}{\exp[-\beta G(x)\{1+\alpha-\alpha G(x)\}]-e^{-\beta}}$$

where $|\alpha|\leq 1, \beta>0, x>0$ and $g(x) = G'(x)$ is the baseline distribution. We refer to this distribution as the Poisson transmuted-G family, in short as $\text{PT-G}(\alpha,\beta)$.

There are many well-known families in the literature. For example, Poisson-G family (Abouelmagd *et al*., 2017), beta-G Poisson family (Gokarna *et al*., 2008), Marshall-Olkin Kumaraswamy-G family (Handique *et al*., 2017), Generalized Marshall-Olkin Kumaraswamy-G family (Chakraborty and Handique 2018), Exponentiated generalized-G Poisson family (Gokarna and Haitham, 2017), beta Kumaraswamy-G family (Handique *et al*., 2017), beta generated Kumaraswamy Marshall-Olkin-G family (Handique and Chakraborty, 2017a), beta generalized Marshall-Olkin Kumaraswamy-G family (Handique and Chakraborty, 2017b), beta generated Marshall-Olkin-Kumaraswamy-G (Chakraborty *et al*., 2018), exponentiated generalized Marshall-Olkin-G family by (Handique *et al*., 2018) and Kumaraswamy generalized Marshall-Olkin-G family (Chakraborty and Handique, 2018) among others.

The primary motivation behind the work is to propose an extended class of distribution that contains the P-G and T-G distribution by adding additional parameter, which covers some important distributions as special and related cases. The PT-G family of distribution appears to be more flexible and could be used for modelling various types of data. For illustration propose we provide pdf of special models of this family in figure 1. It also may work better, in terms of model fitting, than other classes of distributions in certain practical situations.

### 1.1 Important sub models

In this section we provide some special cases of the PT-G family of distributions and list their main distributional characteristics.

> **The PT- exponential (PT-E) distribution**

Let the base line distribution be exponential with parameter $\lambda>0$, $g(x) = \lambda e^{-\lambda x}$ and $G(x) = 1-e^{-\lambda x}$, $x>0$, then for the PT-E model we get the pdf and hrf respectively as

$$f^{\text{PTE}}(x;\alpha,\beta,\lambda) = \frac{\beta \lambda e^{-\lambda x}[1+\alpha-2\alpha(1-e^{-\lambda x})]\exp[-\beta(1-e^{-\lambda x})\{1+\alpha-\alpha(1-e^{-\lambda x})\}]}{1-e^{-\beta}}$$

and $h^{\text{PTE}}(x;\alpha,\beta,\lambda) = \dfrac{\beta\lambda e^{-\lambda x}[1+\alpha-2\alpha(1-e^{-\lambda x})]\exp[-\beta(1-e^{-\lambda x})\{1+\alpha-\alpha(1-e^{-\lambda x})\}]}{\exp[-\beta(1-e^{-\lambda x})\{1+\alpha-\alpha(1-e^{-\lambda x})\}]-e^{-\beta}}$

## ➢ The PT- Weibull (PT-W) distribution

Taking the Weibull distribution (Weibull, 1951) with parameters $\lambda > 0$ and $\theta > 0$ having pdf and cdf $g(x) = \lambda\theta x^{\theta-1}e^{-\lambda x^{\theta}}$ and $G(x) = 1 - e^{-\lambda x^{\theta}}$, $x > 0$ respectively we get the pdf and hrf of PT-W distribution respectively as

$$f^{\text{PTW}}(x;\alpha,\beta,\lambda,\theta) = \dfrac{\beta\lambda\theta x^{\theta-1}e^{-\lambda x^{\theta}}[1+\alpha-2\alpha(1-e^{-\lambda x^{\theta}})]\exp[-\beta(1-e^{-\lambda x^{\theta}})\{1+\alpha-\alpha(1-e^{-\lambda x^{\theta}})\}]}{1-e^{-\beta}}$$

and

$$h^{\text{PTW}}(x;\alpha,\beta,\lambda,\theta) = \dfrac{\beta\lambda\theta x^{\theta-1}e^{-\lambda x^{\theta}}[1+\alpha-2\alpha(1-e^{-\lambda x^{\theta}})]\exp[-\beta(1-e^{-\lambda x^{\theta}})\{1+\alpha-\alpha(1-e^{-\lambda x^{\theta}})\}]}{\exp[-\beta(1-e^{-\lambda x^{\theta}})\{1+\alpha-\alpha(1-e^{-\lambda x^{\theta}})\}]-e^{-\beta}}$$

Obviously a large number of other distributions can be generated by assuming different base line distributions.

### 1.2 Shape of the pdfs

Here we have plotted the pdf of the PT-E and PT-W with chosen parameter values to study the variety of shapes assumed by the family.

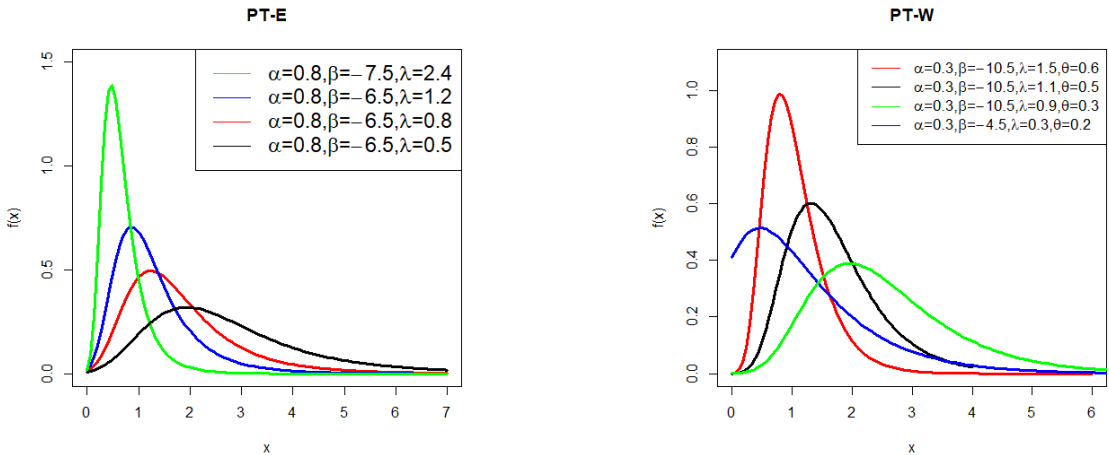

Fig 1: Plots of the pdf of the PT-E and PT-W distributions

## 2. Linear Representation

We express (5) and (6) as infinite series expansion to show that the PT-G can be written as a linear combination of T-G as well as a linear combination of exponentiated-G distributions. These expressions will be helpful to study the mathematical characteristics of the PT-G family.

Using the power series for the exponential function, we can write (5) as

$$f^{PTG}(x;\alpha,\beta) = g^{TG}(x;\alpha)\sum_{i=0}^{\infty}\delta_i [G^{TG}(x;\alpha)]^i \quad (7)$$

$$= \sum_{i=0}^{\infty}\delta'_i \frac{d}{dx}[G^{TG}(x;\alpha)]^{i+1} \quad (8)$$

where $\delta'_i = \frac{(-1)^i \beta^{i+1}}{(1-e^{-\beta})(i+1)i!}$ and $\delta_i = \delta'_i (i+1)$

Using Taylor series expansion cdf of (6) we can write

$$F^{PTG}(x;\alpha,\beta) = \sum_{j=0}^{\infty}\xi_j [G^{TG}(x;\alpha))]^j \quad (9)$$

where $\xi_j = \frac{(-1)^{j+1}\beta^j}{(1-e^{-\beta})j!}$

### 2.1 Moment Generating Function

The moment generating function of PT-G family can be easily expressed in terms of those of the exponentiated T-G distribution using the results of Section 2. For example using equation (8) it can be seen that

$$M_X(s) = E[e^{sX}] = \int_{-\infty}^{\infty} e^{sx} f^{PTG}(x;\alpha,\beta,)\,dx = \int_{-\infty}^{\infty} e^{sx}\sum_{i=0}^{\infty}\delta'_i \frac{d}{dx}[G^{PTG}(x;\alpha)]^{i+1}\,dx$$

$$= \sum_{i=0}^{\infty}\delta'_i \int_{-\infty}^{\infty} e^{sx}\frac{d}{dx}[G^{TG}(x;\alpha)]^{i+1}\,dx = \sum_{i=0}^{\infty}\delta_i M_X(s)$$

where $M_X(s)$ is the mgf of a T-G distribution.

### 2.2 Distribution of Order Statistics

Consider a random sample $X_1, X_2,...,X_n$ from any PT-G distribution. Let $X_{r:n}$ denote the $r^{th}$ order statistic. The pdf of $X_{r:n}$ can be expressed as

$$f_{r:n}(x) = \frac{n!}{(r-1)!(n-r)!} f^{PTG}(x) F^{PTG}(x)^{r-1}\{1-F^{PTG}(x)\}^{n-r}$$

$$= \frac{n!}{(r-1)!(n-r)!} \sum_{m=0}^{n-r} (-1)^m \binom{n-r}{m} f^{PTG}(x)[F^{PTG}(x)]^{m+r-1}$$

The pdf of the $r^{th}$ order statistic for of the PT-G can be derived by using the expansion of its pdf and cdf as

$$f_{r:n}(x) = \frac{n!}{(r-1)!(n-r)!} \sum_{m=0}^{n-r} (-1)^m \binom{n-r}{m} g^{TG}(x;\alpha) \sum_{i=0}^{\infty} \delta_i [G^{TG}(x;\alpha)]^i \left[ \sum_{j=0}^{\infty} \xi_j [G^{TG}(x;\alpha)]^j \right]^{m+r-1}$$

where $\delta_i$ and $\xi_j$ defined above

Using power series raised to power for positive integer $n \ (\geq 1)$ (see Gradshteyn and Ryzhik, 2007)

$$\left( \sum_{i=0}^{\infty} a_i u^i \right)^n = \sum_{i=0}^{\infty} c_{n,i} u^i,$$ where the coefficient $c_{n,i}$ for $i = 1, 2, ...$ are easily obtained from the recurrence equation $c_{n,i} = (i a_0)^{-1} \sum_{m=1}^{i} [m(n+1) - i] a_m c_{n,i-m}$ where $c_{n,0} = a_o^n$.

Now 
$$\left[ \sum_{j=0}^{\infty} \xi_j [G^{TG}(x;\alpha)]^j \right]^{m+r-1} = \sum_{j=0}^{\infty} d_{m+r-1, j} [G^{OMEGP}(x;\beta,\xi)]^j$$

Therefore the density function of the $r^{th}$ order statistics of PT-G distribution can be expressed as

$$f_{r:n}(x) = \frac{n!}{(r-1)!(n-r)!} \sum_{m=0}^{n-r} (-1)^m \binom{n-r}{m} \sum_{i,j=0}^{\infty} \delta_i \, d_{m+r-1, j} [G^{TG}(x;\alpha)]^{i+j} g^{TG}(x;\alpha)$$

$$= \sum_{i,j=0}^{\infty} \psi_{i,j} [G^{TG}(x;\alpha)]^{i+j} g^{TG}(x;\alpha) \tag{10}$$

$$= \sum_{i,j=0}^{\infty} \frac{\psi_{i,j}}{i+j+1} \frac{d}{dx} [G^{TG}(x;\alpha)]^{i+j+1}$$

where $\psi_{i,j} = \frac{n!}{(r-1)!(n-r)!} \sum_{m=0}^{n-r} (-1)^m \binom{n-r}{m} \delta_i \, d_{m+r-1, j}$

## 2.3 Probability Weighted Moments

The probability weighted moments (PWM), first proposed by Greenwood *et al.* (1979), are expectations of certain functions of a random variable whose mean exists. The $(p,q,r)^{th}$ PWM of $T$ is defined by

$$\Gamma_{p,q,r} = \int_{-\infty}^{\infty} x^p F(x)^q [1-F(x)]^r f(x) dx.$$

From equation (7) the $s^{th}$ moment of $T$ can be written as

$$E(X^s) = \int_0^{\infty} x^s f^{PTG}(x;\alpha,\beta) dx = \sum_{i=0}^{\infty} \delta_i \int_0^{\infty} x^s [G^{TG}(x;\alpha)]^i g^{TG}(x;\alpha) dx = \sum_{i=0}^{\infty} \delta_i \Gamma_{s,i,0}$$

where $\delta_i$ is the defined in section 2. Therefore the PWMs of the PT-G can be expressed in terms of linear combination of the PWMs of the T-G distributions.

Proceeding similarly we can express $s^{th}$ moment of the $r^{th}$ order statistic $X_{i:n}$ in a random sample of size $n$ from PT-G on using equation (10) as $E(X^s_{i:n}) = \sum_{i,j=0}^{\infty} \psi_{i,j} \Gamma_{s,i+j,0}$, where $\psi_{i,j}$ defined in above.

### 2.4 Stress-Strength System Reliability

In stress-strength modeling $R = P(X_1 < X_2)$ is a measure of component reliability of the system with random stress $X_1$ and strength $X_2$. It measures the probability that the systems strength $X_2$ is greater than environmental stress $X_1$ applied on that system. The probability of failure of a system is based on the probability of stress exceeding strength, whereas, the reliability of the system is the reversed probability. the system reliability is given by

$$R = P(X_1 < X_2) = P(Stress < Strength) = \int_0^{\infty} f_{Stress}(x) F_{Strength}(x) dx$$

Let $X_1$ and $X_2$ be two independent random variables with $PT-G(x;\alpha_1,\beta_1)$ and $PT-G(x;\alpha_2,\beta_2)$ distributions respectively. Then we have

$$R = \int_0^{\infty} f_1(x;\alpha_1,\beta_1) F_2(x;\alpha_2,\beta_2) dx$$

where $f_1(.)$ and $F_2(.)$ are the pdf and cdf of the PT-G random variables $X_1$ and $X_2$ respectively. Note that the pdf and cdf of $X_1$ and $X_2$ are given by

$$f_1^{PTG}(x;\alpha_1,\beta_1) = g^{TG}(x;\alpha_1) \sum_{i=0}^{\infty} \delta_i^{(1)} [G^{TG}(x;\alpha_1)]^i \text{ and } F^{PTG}(x;\alpha_2,\beta_2) = \sum_{j=0}^{\infty} \xi_j^{(2)} [G^{TG}(x;\alpha_2))]^j$$

Thus $R = \sum_{i=0}^{\infty} \sum_{j=0}^{\infty} \delta_i^{(1)} \xi_j^{(2)} \int_0^{\infty} g^{TG}(x;\alpha_1)[G^{TG}(x;\alpha_1)]^i [G^{TG}(x;\alpha_2))]^j dx$

where $\delta_i^{(1)} = \dfrac{(-1)^i \beta_1^{i+1}}{(1-e^{-\beta_1}) i!}$ and $\xi_j^{(2)} = \dfrac{(-1)^{j+1} \beta_2^j}{(1-e^{-\beta_2}) j!}$

## 2.5 Residual Life and Reversed Residual Life

Let $X$ be a PT-G random variable and $F(x)$ be its cdf (6). The $n^{th}$ moment of the residual life, say $m_n(t) = E[(X-t)^n / X > t]$, $n = 1, 2, \ldots$ uniquely determines $F(x)$. The $n^{th}$ moment of the residual life of $X$ is given by

$$m_n(t) = \frac{1}{1-F(t)} \int_t^\infty (x-t)^n \, dF(x)$$

$$= \frac{1}{1-F(t)} \int_t^\infty \sum_{r=0}^n \binom{n}{r} x^r (-t)^{n-r} f(x) \, dx$$

$$= \frac{1}{1-F(t)} \sum_{i=0}^\infty \delta_i^* \int_t^\infty x^r [G^{TG}(x;\alpha)]^i \, g^{TG}(x;\alpha) \, dx \qquad (11)$$

where $\delta_i^* = \delta_i \sum_{r=0}^n \binom{n}{r} (-t)^{n-r}$.

The $n^{th}$ moment of the reverse residual life, say $M_n(t) = E[(t-X)^n / X \leq t]$, $t > 0$, $n = 1, 2, \ldots$ uniquely determines $F(x)$. We have

$$M_n(t) = \frac{1}{F(t)} \int_0^t (t-x)^n \, dF(x)$$

$$= \frac{1}{F(t)} \int_0^t \sum_{r=0}^n (-1)^r \binom{n}{r} x^r (-t)^{n-r} f(x) \, dx$$

$$= \frac{1}{F(t)} \sum_{i=0}^\infty \delta_i^{**} \int_o^t x^r [G^{TG}(x;\alpha)]^i \, g^{TG}(x;\alpha) \, dx \qquad (12)$$

where $\delta_i^{**} = \delta_i (-1)^n \sum_{r=0}^n \binom{n}{r} t^{n-r}$. The mean residual life (MRL) of $X$ can be obtained by setting $n=1$ in equation (11) and is defined by $m_1(t) = E[(X-t)/X > t]$ also called the life expectation at age $t$ which represents the expected additional life length for a unit which is alive at age $t$. The mean inactivity time (MIT) or mean waiting time (MWT), also called the mean reversed residual life function, is given by $M_1(t) = E[(t-X)/X \leq t]$, $t > 0$ and it represents the waiting time elapsed since the failure of an item on the condition that this failure had occurred in $(0,t)$. The MIT of the PT-G family of distributions can be obtained easily by setting $n=1$ in equation (12).

## 2.6 Rényi entropy

The entropy of a random variable is a measure of uncertainty. The Rényi entropy is defined as

$$I_R(\delta) = (1-\delta)^{-1} \log \left( \int_{-\infty}^\infty f(t)^\delta \, dt \right),$$

where $\delta > 0$ and $\delta \neq 1$. Using the power series for the exponential function, we can write (5) as

$$f^{PTG}(x;\alpha,\beta)^{\delta} = g^{TG}(x;\alpha)^{\delta} \sum_{i=0}^{\infty} \mu_i [G^{TG}(x;\alpha)]^{i\delta}, \text{ where } \mu_i = \frac{(-1)^i \beta^{\delta(i+1)}}{(1-e^{-\beta})^{\delta} i!}$$

Therefore, the Rényi entropy of the PT-G family is given by

Thus
$$I_R(\delta) = (1-\delta)^{-1} \log \left( \int_0^{\infty} g^{TG}(x;\alpha)^{\delta} \sum_{i=0}^{\infty} \mu_i [G^{TG}(x;\alpha)]^{i\delta} dx \right)$$

$$= (1-\delta)^{-1} \log \left( \sum_{i=0}^{\infty} \mu_i \int_0^{\infty} g^{TG}(x;\alpha)^{\delta} [G^{TG}(x;\alpha)]^{i\delta} dx \right)$$

## 2.7 Mean Deviation

Let $X$ be the PT-G random variable with mean $\mu = E(X)$ and median $M = \text{Median}(X) = Q(0.5)$. The mean deviation from the mean $[\delta_{\mu}(X) = E(|X - \mu|)]$ and the mean deviation from the median $[\delta_M(X) = E(|X - M|)]$ can be expressed as

$$\delta_{\mu}(X) = \int_{-\infty}^{\infty} |X - \mu| f(x) dx = \int_{-\infty}^{\mu} (\mu - x) f(x) dx + \int_{\mu}^{\infty} (x - \mu) f(x) dx = 2\mu F(\mu) - 2\Phi(\mu)$$

and
$$\delta_M(X) = \int_{-\infty}^{\infty} |X - M| f(x) dx = \int_{-\infty}^{M} (M - x) f(x) dx + \int_{M}^{\infty} (x - M) f(x) dx = \mu - 2\Phi(M)$$

respectively, where $F(\cdot)$ is the cdf of the PT-G distribution, and $\Phi(t) = \int_{-\infty}^{t} x f(x) dx$

We compute $\Phi(t)$ as follows:

$$\Phi(t) = \sum_{i=0}^{\infty} \delta_i \int_{-\infty}^{t} x g^{TG}(x;\alpha)[G^{TG}(x;\alpha)]^i dx, \text{ where } \delta_i \text{ defined in section 2.}$$

## 3 Maximum Likelihood Estimation

Let $x = (x_1, x_2, ..., x_n)$ be a random sample of size $n$ from PT-G with parameter vector $\boldsymbol{\rho} = (\alpha, \beta, \boldsymbol{\xi})$, where $\boldsymbol{\xi} = (\xi_1, \xi_2, ..., \xi_q)$ is the parameter vector of $G$. The log-likelihood function is written as

$$\ell = \ell(\boldsymbol{\rho}) = n \log \beta - n \log(1 - e^{-\beta}) + \sum_{i=1}^{n} \log[g(x_i, \boldsymbol{\xi})] + \sum_{i=1}^{n} \log[1 + \alpha - 2\alpha G(x_i, \boldsymbol{\xi})]$$

$$- \beta \sum_{i=1}^{n} G(x_i, \boldsymbol{\xi})[1 + \alpha - \alpha G(x_i, \boldsymbol{\xi})]$$

This log-likelihood function cannot be solved analytically because of its complex form but it can be maximized numerically by employing global optimization methods available with the software's R. By taking the partial derivatives of the log-likelihood function with respect to the parameter $\alpha$ and $\beta$, we obtain the components of the score vector $U_\rho = (U_\alpha, U_\beta, U_\xi)$.

The asymptotic variance-covariance matrix of the MLEs of parameters can obtained by inverting the Fisher information matrix $I(\rho)$ which in turn can be derived using the second partial derivatives of the log-likelihood function with respect to each parameter. The $ij^{th}$ elements of $I_n(\rho)$ are given by

$$I_{ij} = -E[\partial^2 l(\rho)/\partial \rho_i \partial \rho_j], \quad i, j = 1, 2+q.$$

The exact evaluation of the above expectations may be cumbersome. In practice one can estimate $I_n(\rho)$ by the observed Fisher's information matrix $\hat{I}_n(\hat{\rho}) = (\hat{I}_{ij})$ defined as

$$\hat{I}_{ij} \approx \left(-\partial^2 l(\rho)/\partial \rho_i \partial \rho_j\right)_{\eta=\hat{\eta}}, \quad i, j = 1, 2+q.$$

Using the general theory of MLEs under some regularity conditions on the parameters as $n \to \infty$ the asymptotic distribution of $\sqrt{n}(\hat{\rho} - \rho)$ is $N_k(0, V_n)$ where $V_n = (v_{jj}) = I_n^{-1}(\rho)$. The asymptotic behaviour remains valid if $V_n$ is replaced by $\hat{V}_n = \hat{I}^{-1}(\hat{\rho})$. Using this result large sample standard errors of $j^{th}$ parameter $\rho_j$ is given by $\sqrt{\hat{v}_{jj}}$.

## 5 Real life applications

Here we consider fitting of two failure time data sets to show that the distributions from the proposed $PT-E(\alpha, \beta, \lambda)$ family can provide better model than the corresponding distributions exponential (Exp), moment exponential (ME), Marshall-Olkin exponential (MO-E) (Marshall and Olkin, 1997), generalized Marshall-Olkin exponential (GMO-E) (Jayakumar and Mathew, 2008), Kumaraswamy exponential (Kw-E) (Cordeiro and de Castro, 2011), Beta exponential (BE) (Eugene *et al.*, 2002), Marshall-Olkin Kumaraswamy exponential (MOKw-E) (Handique *et al.*, 2017), Kumaraswamy Marshall-Olkin exponential (KwMO-E) (Alizadeh *et al.*, 2015), beta Poisson exponential (Handique *et al.*, 2020) and Kumaraswamy Poisson exponential (Chakraborty *et al.*, 2020) distribution.

We have considered some well known model selection criteria namely the AIC, BIC, CAIC and HQIC and the Kolmogorov-Smirnov (K-S) statistics, Anderson-Darling (A) and Cramer von-mises (W) for goodness of fit to compare the fitted models. We have also provided the asymptotic standard errors and confidence intervals of the mles of the parameters for each competing models. Visual comparison fitted density and the fitted cdf are presented in Figures 4 and 5. These plots reveal that the proposed distributions provide a good fit to these data.

The first data is about survival times (in days) of 72 guinea pigs infected with virulent tubercle bacilli, observed and reported by Bjerkedal (1960), while the second one represents the lifetime data relating to relief times (in minutes) of patients receiving an analgesic. The data was reported by Gross and Clark (1975) and it has twenty (20) observations. The descriptive statistics of the data sets are tabulated in Table 1 and both the data are positively skewed.

The total time on test (TTT) plot (see Aarset, 1987) is a technique to extract the information about the shape of the hazard function. A straight diagonal line indicates constant hazard for the data set, where as a convex (concave) shape implies decreasing (increasing) hazard. The TTT plots for the data sets Fig. 3 indicate that the both data sets have increasing hazard rate.

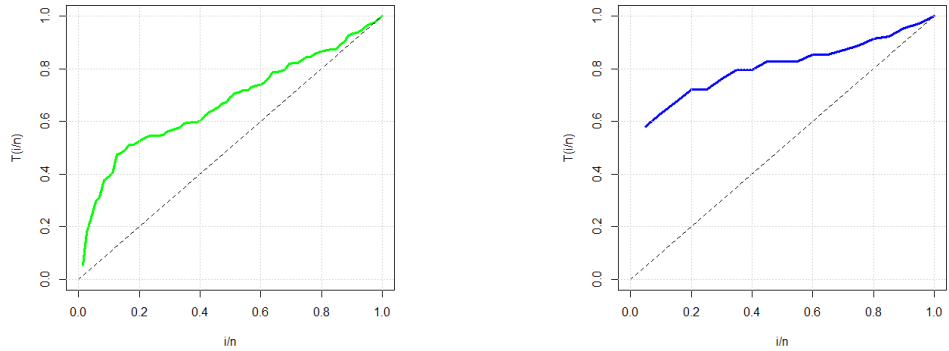

**Fig: 3** TTT-plots for the Data set I and Data set II

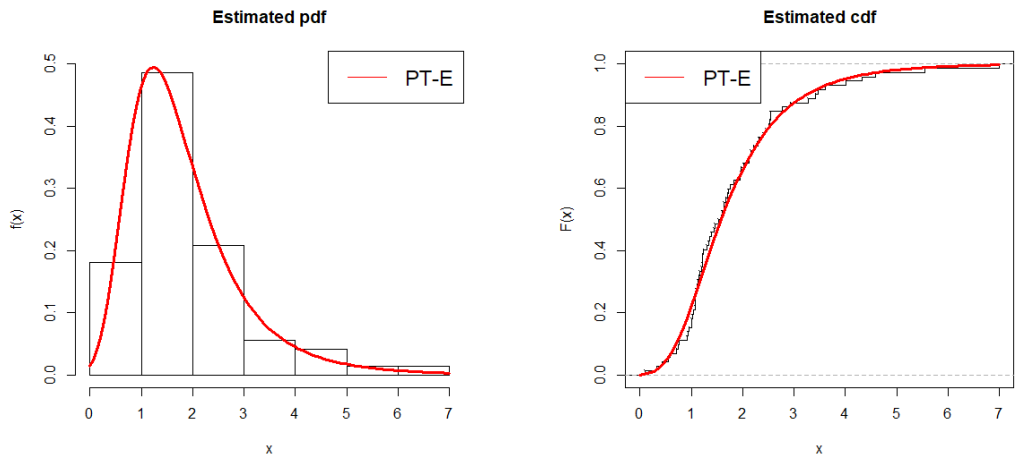

**Fig: 4** Plots of the observed histogram and estimated pdf on left and observed ogive and estimated cdf for the PT-E model for data set I

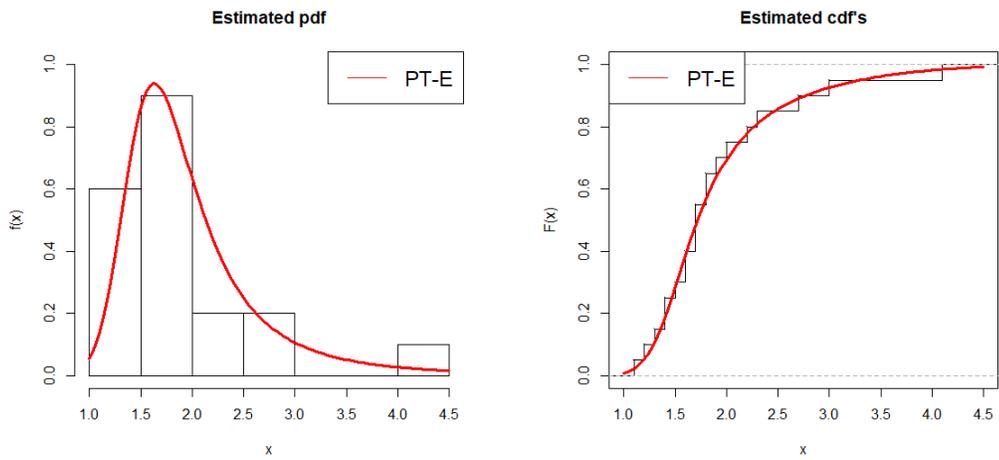

**Fig: 5** Plots of the observed histogram and estimated pdf on left and observed ogive and estimated cdf for the PT-E model for data set II

**Table 1:** Descriptive Statistics for the data sets I and II

| Data Sets | n | Min. | Mean | Median | s.d. | Skewness | Kurtosis | 1st Qu. | 3rd Qu. | Max. |
|---|---|---|---|---|---|---|---|---|---|---|
| I | 72 | 0.100 | 1.851 | 1.560 | 1.200 | 1.788 | 4.157 | 1.080 | 2.303 | 7.000 |
| II | 20 | 1.100 | 1.900 | 1.700 | 0.704 | 1.592 | 2.346 | 1.475 | 2.050 | 4.100 |

**Table 2 (a):** MLEs, standard errors, confidence intervals (in parentheses) values for the guinea pigs survival time's data set I

| Models | $\hat{a}$ | $\hat{b}$ | $\hat{\alpha}$ | $\hat{\beta}$ | $\hat{\lambda}$ |
|---|---|---|---|---|---|
| Exp $(\lambda)$ | --- | --- | --- | --- | 0.540 (0.063) (0.42, 0.66) |
| ME $(\lambda)$ | --- | --- | --- | --- | 0.925 (0.077) (0.62, 1.08) |
| MO-E $(\alpha, \lambda)$ | --- | --- | 8.778 (3.555) (1.81, 15.74) | --- | 1.379 (0.193) (1.00, 1.75) |
| GMO-E $(b, \alpha, \lambda)$ | --- | 0.179 (0.070) (0.04, 0.32) | 47.635 (44.901) (0, 135.64) | --- | 4.465 (1.327) (1.86, 7.07) |
| Kw-E $(a, b, \lambda)$ | 3.304 (1.106) (1.13, 5.47) | 1.100 (0.764) (0, 2.59) | --- | --- | 1.037 (0.614) (0, 2.24) |
| B-E $(a, b, \lambda)$ | 0.807 (0.696) (0, 2.17) | 3.461 (1.003) (1.49, 5.42) | --- | --- | 1.331 (0.855) (0, 3.01) |
| MOKw-E $(\alpha, a, b, \lambda)$ | 2.716 (1.316) (0.14, 5.29) | 1.986 (0.784) (0.449, 3.52) | 0.008 (0.002) (0.004, 0.01) | --- | 0.099 (0.048) (0, 0.19) |
| KwMO-E $(a, b, \alpha, \lambda)$ | 3.478 (0.861) (1.79, 5.17) | 3.306 (0.779) (1.78, 4.83) | 0.373 (0.136) (0.11, 0.64) | --- | 0.299 (1.112) (0, 2.48) |
| BP-E $(a, b, \beta, \lambda)$ | 3.595 (1.031) (1.57, 5.62) | 0.724 (1.590) (0, 3.84) | --- | 0.014 (0.010) (0, 0.03) | 1.482 (0.516) (0.47, 2.49) |
| KwP-E $(a, b, \beta, \lambda)$ | 3.265 (0.991) (1.32, 5.21) | 2.658 (1.984) (0, 6.55) | --- | 4.001 (5.670) (0, 15.11) | 0.177 (0.226) (0, 0.62) |
| PT-E $(\alpha, \beta, \lambda)$ | --- | --- | 0.813 (0.182) (0.45, 1.17) | -6.587 (1.448) (0, -3.74) | 0.841 (0.192) (0.46, 1.22) |

**Table 2 (b):** Log-likelihood, AIC, BIC, CAIC, HQIC, A, W and KS (*p*-value) values for the guinea pigs survival times data set I

| Models | AIC | BIC | CAIC | HQIC | A | W | KS (*p*-value) |
|---|---|---|---|---|---|---|---|
| Exp $(\lambda)$ | 234.63 | 236.91 | 234.68 | 235.54 | 6.53 | 1.25 | 0.27 (0.06) |
| ME $(\lambda)$ | 210.40 | 212.68 | 210.45 | 211.30 | 1.52 | 0.25 | 0.14 (0.13) |
| MO-E $(\alpha, \lambda)$ | 210.36 | 214.92 | 210.53 | 212.16 | 1.18 | 0.17 | 0.10 (0.43) |
| GMO-E $(b, \alpha, \lambda)$ | 210.54 | 217.38 | 210.89 | 213.24 | 1.02 | 0.16 | 0.09 (0.51) |
| Kw-E $(a, b, \lambda)$ | 209.42 | 216.24 | 209.77 | 212.12 | 0.74 | 0.11 | 0.08 (0.50) |
| B-E $(a, b, \lambda)$ | 207.38 | 214.22 | 207.73 | 210.08 | 0.98 | 0.15 | 0.11 (0.34) |
| MOKw-E $(\alpha, a, b, \lambda)$ | 209.44 | 218.56 | 210.04 | 213.04 | 0.79 | 0.12 | 0.10 (0.44) |
| KwMO-E $(a, b, \alpha, \lambda)$ | 207.82 | 216.94 | 208.42 | 211.42 | 0.61 | 0.11 | 0.08 (0.73) |
| BP-E $(a, b, \beta, \lambda)$ | 205.42 | 214.50 | 206.02 | 209.02 | 0.55 | 0.08 | 0.09 (0.81) |
| KwP-E $(a, b, \beta, \lambda)$ | 206.63 | 215.74 | 207.23 | 210.26 | 0.48 | 0.07 | 0.09 (0.79) |
| **PT-E $(\alpha, \beta, \lambda)$** | **202.09** | **208.92** | **202.44** | **204.81** | **0.36** | **0.05** | 0.07 (**0.86**) |

**Table 3 (a):** MLEs, standard errors, confidence intervals (in parentheses) values for the relief times of patients receiving an analgesic failure time data set II

| Models | $\hat{a}$ | $\hat{b}$ | $\hat{\alpha}$ | $\hat{\beta}$ | $\hat{\lambda}$ |
|---|---|---|---|---|---|
| Exp $(\lambda)$ | --- | --- | --- | --- | 0.526 (0.117) (0.29, 0.75) |
| ME $(\lambda)$ | --- | --- | --- | --- | 0.950 (0.150) (0.66, 1.24) |
| MO-E $(\alpha, \lambda)$ | --- | --- | 54.474 (35.582) (0, 124.21) | --- | 2.316 (0.374) (1.58, 3.04) |
| GMO-E $(b, \alpha, \lambda)$ | --- | 0.519 (0.256) (0.02, 1.02) | 89.462 (66.278) (0, 219.37) | --- | 3.169 (0.772) (1.66, 4.68) |
| Kw-E $(a, b, \lambda)$ | 83.756 (42.361) (0.73, 166.78) | 0.568 (0.326) (0, 1.21) | --- | --- | 3.330 (1.188) (1.00, 5.66) |
| B-E $(a, b, \lambda)$ | 81.633 (120.41) (0, 317.63) | 0.542 (0.327) (0, 1.18) | --- | --- | 3.514 (1.410) (0.75, 6.28) |
| MOKw-E $(\alpha, a, b, \lambda)$ | 33.232 (57.837) (0, 146.59) | 0.571 (0.721) (0, 1.98) | 0.133 (0.332) (0, 0.78) | --- | 1.669 (1.814) (0, 5.22) |
| KwMO-E $(a, b, \alpha, \lambda)$ | 34.826 (22.312) (0, 78.56) | 0.299 (0.239) (0, 0.76) | 28.868 (9.146) (10.94, 46.79) | --- | 4.899 (3.176) (0, 11.12) |
| BP-E $(a, b, \beta, \lambda)$ | 13.396 (1.494) (10.46, 16.32) | 9.600 (1.091) (7.46, 11.73) | --- | 1.965 (0.341) (1.29, 2.63) | 0.244 (0.037) (0.17, 0.32) |
| KwP-E $(a, b, \beta, \lambda)$ | 11.837 (6.493) (0, 24.56) | 3.596 (2.392) (0, 8.28) | --- | 5.983 (1.470) (3.10, 8.86) | 0.225 (0.098) (0.03, 0.42) |
| PT-E $(\alpha, \beta, \lambda)$ | --- | --- | 0.301 (0.037) (0.22, 0.37) | -9.997 (3.336) (0, -3.54) | 1.555 (0.241) (1.08, 2.02) |

**Table 3 (b):** Log-likelihood, AIC, BIC, CAIC, HQIC, A, W and KS (*p*-value) values for the relief times of patients receiving an analgesic failure time data set II

| Models | AIC | BIC | CAIC | HQIC | A | W | KS (*p*-value) |
|---|---|---|---|---|---|---|---|
| Exp $(\lambda)$ | 67.67 | 68.67 | 67.89 | 67.87 | 4.60 | 0.96 | 0.44 (0.004) |
| ME $(\lambda)$ | 54.32 | 55.31 | 54.54 | 54.50 | 2.76 | 0.53 | 0.32 (0.07) |
| MO-E $(\alpha, \lambda)$ | 43.51 | 45.51 | 44.22 | 43.90 | 0.81 | 0.14 | 0.18 (0.55) |
| GMO-E $(b, \alpha, \lambda)$ | 42.75 | 45.74 | 44.25 | 43.34 | 0.51 | 0.08 | 0.15 (0.78) |
| Kw-E $(a, b, \lambda)$ | 41.78 | 44.75 | 43.28 | 42.32 | 0.45 | 0.07 | 0.14 (0.86) |
| B-E $(a, b, \lambda)$ | 43.48 | 46.45 | 44.98 | 44.02 | 0.70 | 0.12 | 0.16 (0.80) |
| MOKw-E $(\alpha, a, b, \lambda)$ | 41.58 | 45.54 | 44.25 | 42.30 | 0.60 | 0.11 | 0.14 (0.87) |
| KwMO-E $(a, b, \alpha, \lambda)$ | 42.88 | 46.84 | 45.55 | 43.60 | 1.08 | 0.19 | 0.15 (0.86) |
| BP-E $(a, b, \beta, \lambda)$ | 38.07 | 42.02 | 40.73 | 38.78 | 0.39 | 0.06 | 0.14 (0.91) |
| KwP-E $(a, b, \beta, \lambda)$ | 38.32 | 42.28 | 40.98 | 39.04 | 0.41 | 0.05 | 0.13 (0.93) |
| **PT-E $(\alpha, \beta, \lambda)$** | **36.84** | **39.81** | **38.34** | **37.38** | **0.37** | **0.04** | **0.11 (0.95)** |

The mle's of the parameters with corresponding standard errors in the parentheses for all the fitted models along are given in Table 2(a) and Table 3(a) for data set I and data set II respectively. While the various model selection criteria namely the AIC, BIC, CAIC, HQIC, A, W and KS statistic with p-value for the fitted models of the data sets I and II are presented respectively in Table 2(b) and Table 3(b). From these findings based on the lowest values different criteria the PT-E is found to be a better model than the models Exp, ME, MO-E, GMO-E, Kw-E, B-E, MOKw-E, KwMO-E, BP-E and KwP-E for both the data sets. A visual comparison of the closeness of the fitted density with the observed histogram and fitted cdf with the observed ogive of the data sets I and II are presented in the Figures 4 and 5 respectively also indicate that the proposed distributions provide comparatively closer fit to these data sets.

## 6 Conclusion

We present a new Poisson transmuted-G (PT-G) family of distributions, which extends the transmuted family by adding one extra shape parameter. Many well-known distributions emerge as special cases of the PT-G family by using special parameter values. The mathematical properties of the new family including explicit expansions for the moments generating function, order statistics, Probability weighted moments, stress-strength reliability, residual life, reversed residual life, Rényi entropy and mean deviation are provided. The model parameters are estimated by the maximum likelihood method. It is shown, by means of two real data sets that special cases of the PT-G family can give better fits than other models generated by well-known families.

**References**

1. Aarset MV (1987) How to identify a bathtub hazard rate. IEEE Transactions on Reliability 36: 106-108.

2. Abouelmagd THM, Hamed MS, Ebraheim AN (2017) The Poisson-G family of distributions with applications. Pak.j.stat.oper.res. XIII: 313-326.

3. Alzaatreh A, Lee C, Famoye F (2013) A new method for generating families of continuous distributions. Metron. 71: 63-79.

4. Alizadeh M, Tahir MH, Cordeiro GM, Zubai M, Hamedani GG (2015) The Kumaraswamy Marshal-Olkin family of distributions. Journal of the Egyptian Mathematical Society, 23: 546-557.

5. Bjerkedal T (1960) Acquisition of resistance in Guinea pigs infected with different doses of virulent tubercle bacilli. American Journal of Hygiene 72: 130-148.

6. Chakraborty S, Handique L (2017) The generalized Marshall-Olkin-Kumaraswamy-G family of distributions. Journal of Data Science 15: 391-422.

7. Chakraborty S, Handique L (2018) Properties and data modelling applications of the Kumaraswamy generalized Marshall-Olkin-G family of distributions. Journal of Data Science 16: 605-620.

8. Chakraborty S, Handique L, Jamal F (2020) The Kumaraswamy Poisson-G family of distribution: its properties and applications. Annals of Data Science (published online).

9. Chakraborty S, Handique L, Ali MM (2018) A new family which integrates beta Marshall-Olkin-G and Marshall-Olkin-Kumaraswamy-G families of distributions. Journal of Probability and Statistical Science 16: 81-101.

10. Cordeiro GM, De Castro M (2011) A new family of generalized distributions. J. Stat. Comput. Simul. 81: 883-893.

11. Eugene N, Lee C, Famoye F (2002) Beta-normal distribution and its applications. Commun Statist Theor Meth, 31: 497-512.

12. Gokarna RA, Sher BC, Hongwei L, Alfred AA (2018) On the Beta-G Poisson family. Annals of Data Science. https://doi.org/10.1007/s40745-018-0176-x.

13. Gokarna RA, Haitham MY (2017) The exponentiated generalized-G Poisson family of distributions. Stochastics and Quality Control 32:7-23.

14. Greenwood JA, Landwehr JM, Matalas NC, Wallis JR (1979) Probability weighted moments: definition and relation to parameters of several distributions expressible in inverse form. Water Resour Res, 15: 1049-1054.

15. Gradshteyn, I.S., Ryzhik, I.M., (2007) Tables of Integrals, Series, and Products, AcademicPress, New York.

16. Gross AJV, Clark A (1975) Survival distributions: Reliability applications in the biometrical sciences. John Wiley, New York, USA.

17. Handique L, Chakraborty S (2017a) A new beta generated Kumaraswamy Marshall-Olkin-G family of distributions with applications. Malaysian Journal of Science 36: 157-174.

18. Handique L, Chakraborty S (2017b) The Beta generalized Marshall-Olkin Kumaraswamy-G family of distributions with applications. Int. J. Agricult. Stat. Sci. 13: 721-733.

19. Handique L, Chakraborty S, Hamedani GG (2017) The Marshall-Olkin-Kumaraswamy-G family of distributions. Journal of Statistical Theory and Applications, 16: 427-447.

20. Handique L, Chakraborty S, Ali MM (2017) The Beta-Kumaraswamy-G family of distributions. Pakistan Journal of Statistics, 33: 467-490.

21. Handique L, Chakraborty S, Thiago AN (2018) The exponentiated generalized Marshall-Olkin family of distributions: Its properties and applications. Annals of Data Science 1-21. https://doi.org/10.1007/s40745-018-0166-z.

22. Handique L, Chakraborty S, Jamal F (2019) Beta Poisson-G family of distribution: Its properties and application with failure time data. Thailand Statistician (Submitted).

23. Haq MA, Handique L, Chakraborty S (2018) The odd moment exponential family of distributions: Its properties and applications. International Journal of Applied Mathematics and Statistics 57: 47-62.

24. Marshall A, Olkin I (1997) A new method for adding a parameter to a family of distributions with applications to the exponential and Weibull families. Biometrika, 84: 641-652.

25. Jayakumar K, Mathew, T (2008) On a generalization to Marshall-Olkin scheme and its application to Burr type XII distribution. Stat Pap, 49: 421-439.

26. Weibull, W (1951) A statistical distribution functions of wide applicability. J. Appl. Mech.-Trans. ASME, 18: 293-297.